\documentclass[12pt,a4paper]{article}
\usepackage[intlimits]{amsmath}
\usepackage{amsfonts,amssymb,amscd,amsthm}
\usepackage{cite}
\input xy
\xyoption{all}

\def\ind{\operatorname{ind}}
\def\Ell{\operatorname{Ell}}

\DeclareMathOperator{\Con}{Con}
\def\RR{\mathbb{R}}
\def\lra{\longrightarrow}
\def\a{\alpha}
\def\be{\beta}
\def\ga{\gamma}

\def\e{\varepsilon}
\def\si{\sigma}
\def\ph{\varphi}
\def\om{\omega}
\def\Om{\Omega}
\def\La{\Lambda}
\DeclareMathOperator\im{Im}
\DeclareMathOperator\supp{supp}
\DeclareMathOperator\const{const}
\DeclareMathOperator\id{id}
\newcommand{\op}[1]{\operatorname{#1}}
\def\cK{{\mathcal{K}}}
\def\cM{{\mathcal{M}}}
\def\ov{\overline}
\def\ovs{\overset}
\def\wt{\widetilde}

\def\pa{\partial}
\newcommand\pd[2]{\frac{\pa #1}{\pa #2}}
\let\rom\textup

\newcommand{\BL}{\biggl}
\newcommand{\BR}{\biggr}
\newcommand{\Vect}{{\rm Vect}}
\let\dokaend\qed
\newcommand{\abs}[1]{\lvert#1\rvert}

\newcommand{\norm}[1]{\left\| #1 \right\|}
\def\Si{\Sigma}
\newtheorem{corollary}{Corollary}[section]
\newtheorem{theorem}{Theorem}[section]
\newtheorem{lemma}{Lemma}[section]
\newtheorem{proposition}{Proposition}[section]

\theoremstyle{definition}

\newtheorem{remark}{Remark}[section]

\newtheorem{definition}{Definition}[section]

\title{\bf On the homotopy classification of elliptic operators
on stratified manifolds\thanks{Supported in part by RFBR grants 05-01-00982
and 06-01-00098, by President of the Russian Federation grant
MK-1713.2005.1, and by the DFG project 436 RUS
113/849/0-1\protect\raisebox{1pt}{\circledR}\ ``K-theory and Noncommutative
Geometry of Stratified Manifolds."}}

\author{{\bf V.E.~Nazaikinskii, A.Yu.~Savin, and B.Yu.~Sternin}}

\date{August 14, 2006}
\begin{document}
\maketitle

\begin{abstract}
We find the stable homotopy classification
of elliptic operators on stratified manifolds. Namely, we
establish an isomorphism of the set of elliptic
operators modulo stable homotopy and the $K$-homology group of the singular
manifold. As a corollary, we obtain an explicit
formula for the obstruction of Atiyah--Bott type to making
interior elliptic operators Fredholm.
\end{abstract}

\tableofcontents

\section{Introduction}

In the classical paper~\cite{Ati4}, Atiyah
observed that the elements of $K$-homology groups of a compact space
$X$ are realized as \emph{abstract elliptic operators} on $X$.
(Recall that the latter are Fredholm operators acting in
$C(X)$-modules and commuting with
multiplications by functions modulo compact operators.)
Moreover, Kasparov~\cite{Kas1} and Brown--Douglas--Fillmore~\cite{BDF1}
showed that one can obtain not only elements
of $K$-homology groups but also $K$-homology as a generalized
homology theory if one factorizes abstract
elliptic operators by a special equivalence relation known as stable
homotopy.

However, it turns out that, at least for sufficiently smooth spaces, the
$K$-homology group can already be obtained if, instead of
abstract elliptic operators, we restrict ourselves to \emph{differential}
or \emph{pseudodifferential} operators, which are natural in
the theory of partial differential equations. For
example, on closed smooth manifolds it suffices to consider
pseudodifferential operators ({$\psi$DO} in what follows), and if the
manifold is equipped with a $spin^c$-structure, then
one can take the class of (twisted) Dirac operators. This example serves as
a motivation for the natural problem of comparing the $K$-homology group
with the group generated by elliptic operators on singular
manifolds (cf.\ Singer's problem in~\cite{Sin2}).

The main theorem of the present paper
establishes a group isomorphism
\begin{equation}\label{mina}
\Ell(X) \simeq K_0(X)
\end{equation}
on any compact stratified manifold $X$ with arbitrary number of
strata, where $\Ell(X)$ is the group generated by elliptic
pseudodifferential operators on $X$ modulo stable homotopy. Special cases
of such isomorphisms were obtained
in~\cite{SaSt10,Sav8,NSScS14,DeLe1,MeRo2} for manifolds with two strata.

The isomorphism~\eqref{mina} enables one to apply topological methods
of $K$-homology theory in elliptic theory. As examples of
such applications, in Sec.~\ref{vosem} we
compute the obstruction of Atiyah--Bott type to making interior elliptic
operators Fredholm and obtain a generalization of
the cobordism invariance of the index.

Apart from already mentioned applications to elliptic
operators, the isomorphism~ \eqref{mina} has an interesting
interpretation in the framework of noncommutative geometry. More precisely,
the groups in \eqref{mina} can be represented as the left- and
right-hand sides in the Baum--Connes
conjecture~\cite{BaCo1}. The point is that from the viewpoint of
noncommutative geometry the algebra of {$\psi$DO} on a stratified manifold
is related to a certain groupoid (see~\cite{LaNi1,Nis99}). Moreover, the
group $\Ell(X)$ is related to the $K$-group of the groupoid
$C^*$-algebra~\cite{LMN1}. The Baum--Connes conjecture claims that the
latter $K$-group is isomorphic to the topological $K$-group of the
classifying space of the groupoid (see~\cite{Tu1}).
Explicit computations for simplest stratified manifolds show that
the $K$-group of the classifying space is isomorphic to $K_0(X)$, that
is, the right-hand side of~\eqref{mina}. It
would be of interest to make further
comparison of~\eqref{mina} with the Baum--Connes map.

We are grateful to Prof.~T.~Fack (Lion) for pointing out the relation with
the Baum--Connes isomorphism for groupoids. The results of this paper were
presented at the International conference ``Workshop
on index theory'' in M\"unster (November 4--8,
2005).

\section{Stratified manifolds and {$\psi$DO}}

In this section, we define a class of manifolds and a class of
pseudodifferential operators on them and study ellipticity and Fredholm
property for such operators. These topics are well known in the literature
(e.g., see~\cite{PlSe1,PlSe6,CMSch1}). Hence the exposition will be rather
concise.

\subsection{Stratified manifolds}

In the sequel, we adopt the following terminology. A \textit{manifold with singularities}
is a triple
$$
\pi: M\longrightarrow \mathcal{M},
$$
where $\cM$ is a Hausdorff space, $M$ is a manifold with
corners,\footnote{Recall that a manifold of dimension $n$ with corners is a
Hausdorff space locally homeomorphic to the product
$\overline{\mathbb{R}}^{k}_+\times \mathbb{R}^{n-k}$, $0 \le k\le n$.}
and $\pi$ is a continuous projection. The manifold
$M$ is called the \emph{blow-up} of $\cM$. We shall not discuss the
uniqueness of the blow-up and always assume when speaking about
manifolds with singularities that the triple $(\pi:M\lra\cM)$ is fixed.

A \textit{diffeomorphism} of such manifolds is a pair $f:\mathcal{M}_1\to
\mathcal{M}_2$, $\widetilde{f}:M_1\to M_2$ such that
$f$ is a homeomorphism, $\wt f$ is a diffeomorphism, and the
diagram
\begin{equation*}
\begin{CD}
M_1@>\wt f>>M_2\\
@V\pi_1VV @VV\pi_2V\\
\cM_1@>>f>\cM_2,
\end{CD}
\end{equation*}
where $\pi_{1,2}$ are the natural projections, commutes. A special class of manifolds with
singularities can be obtained as follows: given a smooth Riemannian metric
on a manifold $M$ with corners, nondegenerate in the interior and
possibly degenerate at the boundary, define $\cM$
as the quotient of $M$ by the following equivalence
relation determined by the metric: two points are
equivalent if the distance between them in our metric is zero.
(Of course, the Hausdorff property of $\cM$ should be additionally assumed.)

Now let us describe the class of manifolds
with singularities
to be studied in the present paper, namely,
stratified manifolds. The definition is by induction.

\begin{definition}\label{def0}
A \textit{filtration of length $k$} on a topological space $\cM$ is a
sequence
\begin{equation}\label{filter}
\mathcal{M}=\mathcal{M}_k\supset \mathcal{M}_{k-1}\supset\ldots\supset
\mathcal{M}_{0}
\end{equation}
of embedded subspaces, closed in $\cM$, such that each $\cM_j$ is
contained in the closure of $\cM_{j+1}^\circ=\cM_{j+1}\setminus\cM_j$,
$j=0,\dotsc,k-1$.
\end{definition}

\begin{definition}[base of induction]\label{def1}
\textit{A stratified manifold of length} $0$ is an arbitrary smooth manifold.

In this case, we set $M=\cM$ and $\pi=\id$, and the
blow-up $M$ is a smooth manifold without boundary.
\end{definition}

\begin{definition}[inductive step]\label{odintri}
A \textit{stratified manifold of length} $k>0$ is a Hausdorff space
$\cM$ equipped with a filtration~\eqref{filter}
such that
\begin{enumerate}
\item[1.] $\cM_0$ is a smooth manifold.
\item[2.]
$\cM\setminus\cM_0$ is equipped with the structure
of a stratified manifold of length $k-1$
with respect to the filtration
\begin{equation*}
\mathcal{M}\setminus\cM_0=\mathcal{M}_r\setminus\cM_0\supset
\mathcal{M}_{k-1}\setminus\cM_0\supset\ldots\supset
\mathcal{M}_{1}\setminus\cM_0.
\end{equation*}
\item[3.] Over $\cM_0$ there is a given bundle with fiber $K_\Om$, where $K_\Om$ is a cone with base
$\Om$ that is a compact stratified manifold of
length less than or equal to $k-1$.
There is also a given homeomorphism of a
neighborhood $U\subset \cM$ of $\cM_0$ onto a
neighborhood of the subbundle formed by the cone tips; this
homeomorphism restricts to the identity map on
$\cM_0$.
\item[4.] The structure described in Condition
3 should be compatible with that in Condition 2 on
$\cM\setminus\cM_0$ in the sense described below.
\end{enumerate}
\end{definition}

First, note that one can prove by induction that the definition implies that
\begin{itemize}
\item The sets $\mathcal{M}_j\setminus \mathcal{M}_{j-1}\simeq
M_j^\circ$ (called open strata) are smooth manifolds for all
$j=0,\dotsc,k$.\footnote{We set $\mathcal{M}_{-1}=\emptyset$ for
convenience. These sets are the interior parts of the corresponding
blow-ups, which are smooth manifolds with corners.} \item Given
any open stratum $M_j^\circ$, $j<k$, there exists a neighborhood
$U_j\subset \mathcal{M}$ homeomorphic to a bundle with fiber
$K_{\La_j}$, where the base of the cone $\La_j$ is a stratified manifold
of length $\le k-j-1$.
\end{itemize}

Let us give the precise statement of Conditions 3 and 4.

The cone $K_\Om$ in Condition 3 is the set
\begin{equation*}
K_\Om=\{\ov\RR_+\times\Om\}\bigm\slash\{\{0\}\times\Om\}.
\end{equation*}
We require that the transition functions of the bundle with fiber
$K_\Omega$ leave the variable $r\in \overline{\mathbb{R}}_+$ invariant and
are induced by \emph{diffeomorphisms of the stratified manifold} $\Omega$
of length $\le k-1$, i.e., diffeomorphisms of manifolds with singularities
preserving the stratification and the fibrations of neighborhoods of strata
into cones. Hence we actually require that our bundle is
obtained from a bundle over $\cM_0$ with fiber $\Om$ by
fiberwise conification.

We claim that the complement $U\setminus\cM_0$ is
equipped with a natural structure of a stratified manifold of length
$\le k-1$ as an open subset of the bundle over $\cM_0$
with fiber
\begin{equation*}
K_\Om^\circ\simeq\RR_+\times\Om
\end{equation*}
equal to the cone minus the cone tip. Indeed, the local trivializations of
this bundle are $V\times\RR_+\times\Om$, where $V\subset\RR^l$, $l=\dim
\mathcal{M}_0$ is a coordinate patch on $\cM_0$. Therefore, they are
stratified manifolds of the same length as $\Om$, i.e., $\le k-1$.
Now we can take $V\times\RR_+\times\Om_j$ for strata in
$V\times\RR_+\times\Om$, where $\Om_j$ are the corresponding strata in
$\Om$, and bundles with conical
fibers of neighborhoods of these strata are obtained from the
bundles of neighborhoods of the corresponding strata
in $\Om$ by taking products with $V\times\RR_+$.

\medskip

Let us clarify Condition 4. By virtue of the preceding,
$U\setminus\cM_0$ is equipped with two stratified
manifold structures: the first is the restriction of the corresponding
structure to $\cM\setminus\cM_0$, and the
second comes from the bundle. \textsl{The
compatibility condition requires that these two structures coincide} (i.e.,
the identity map is a diffeomorphism).

\bigskip

Finally, let us define the blow-up of $\cM$. Let $\wt\pi:\wt
M\lra\cM\setminus\cM_0$ be the projection of the blow-up
$\cM\setminus\cM_0$ (it is defined by the inductive hypothesis). The
blow-up $M$ of $\cM$ is obtained if we augment
$\wt M$ with a certain set ``at infinity over $\cM_0$.''
Hence to describe $M$ and the projection $\pi:M\lra \cM$ it suffices to study what happens over $U$. We can assume that
$U$ itself is fibered over $\mathcal{M}_0$; then
$\wt\pi^{-1}(U\setminus\cM_0)$ is fibered over $\cM_0$ by
the composition map
\begin{equation*}
\wt\pi^{-1}(U\setminus\cM_0)\ovs{\wt\pi}\lra U\setminus\cM_0
\lra\cM_0.
\end{equation*}
The structure of this bundle is easy to describe in
local trivializations; the bundle
has the form
\begin{equation*}
V\times\RR_+\times\wt\Om\xrightarrow{\id\times\id\times p}
V\times\RR_+\times\Om\longrightarrow V,
\end{equation*}
where $p:\wt\Om\lra\Om$ is the projection of the blow-up of the
cone base. In this local trivialization, the blow-up
$M$ of $\cM$ is simply defined by adding the point $r=0$
to the second factor, i.e., by passing from $\RR_+$ to $\ov\RR_+$;
the projection takes each point $(v,0,\om)$ to $v$.

Now the meaning of Definition~\ref{odintri} is completely
clear.

Each of the strata $\cM_j$ is a stratified manifold itself (of length $j$).
It has a blow-up $M_j$ and the corresponding projection $p_j:M_j\lra
\cM_j$.

\begin{remark}
By Condition 3, conical bundles are
defined over the open strata $M_j^\circ$ in $\cM$. However,
Condition 4 readily implies that these
bundles can be canonically extended to the closed
strata $M_j$. This remark will be important below, since the
operator-valued symbols of our pseudodifferential operators will be defined
on the closed strata.
\end{remark}

\paragraph{Metrics, measures, and $L^2$-spaces.}

We introduce some natural metrics and measures on stratified manifolds.
They are used in the definition of spaces of square integrable
functions.

First, we give an inductive definition of metrics. On a stratified manifold
of length 0 we consider an arbitrary Riemannian metric. Let us
describe the inductive step. Given a metric $d\wt\om^2$ on $\wt\Om$
supplied by the induction hypothesis, we can define a metric locally
on $V\times\RR_+\times\wt\Om$ by the formula
\begin{equation}\label{form}
ds^2=dv^2+dr^2+r^2d\wt\om^2.
\end{equation}
Globally on $M$, the metric is obtained with
the use of a partition of unity from the local
expressions described above in a neighborhood of
$\cM_0$ and the metric $\wt d\rho^2$ defined on
$\cM\setminus\cM_0$ by the induction hypothesis outside a slightly
smaller neighborhood of $\cM_0$. Metrics of this form are called
\textit{edge-degenerate}.

The metric naturally produces the measure
defined as the volume element equal to unity on
an orthonormal
frame. In terms of the inductive formula~\eqref{form}, the corresponding
formula for the volume form is
\begin{equation*}
d\op{vol}=r^ndv\,dr\,d\op{vol}_\Om,
\end{equation*}
where $d\op{vol}_\Om$ is the volume form on $\Om$ (known by the induction
hypothesis) and $n=\dim\Om$ is the dimension of $\Om$.

From now on, all operators on $\cM$ are considered
in the space
$$
L^2(\cM)\equiv L^2(\cM,d\op{vol}),
$$
and the operators on the cone $K_\Om$ are
considered in the space
$$
L^2(K_\Om)\equiv L^2(K_\Om,r^n\, dr\,d\op{vol}_\Om).
$$

\paragraph{The cotangent bundle.}

We define a space $\Vect_\cM$ of vector fields on $\cM$. If
$\cM$ is smooth, then $\Vect_\cM$ is the space of all vector fields on
$\cM$. In the general case, $\Vect_\cM$ is defined as a
$C^\infty(M)$-module locally as follows. Over the product
$V\times\RR_+\times\wt\Om$, the space $\Vect_\cM$ is formed by vector
fields of the form
\begin{equation*}
\theta=a\pd{}v+b\pd{}r+\frac1r\theta_1,
\end{equation*}
where $a$ and $b$ are smooth functions and $\theta_1\in\Vect_\Om$.

The metric $ds^2$ defines a $C^\infty(M)$-valued pairing on $\Vect_\cM$,
and the formula
\begin{equation*}
\langle\ph(\theta),\mu\rangle=ds^2(\theta,\mu)
\end{equation*}
gives a bijection $\varphi$ of the space $\Vect_\cM$ onto some
$C^\infty(M)$-module $\Lambda^1(\cM)\subset\Lambda^1(M)$ of differential
forms on the blow-up $M$. (To see this, it suffices to note that the
embedding $\Vect_M\subset\Vect_\cM$ is dense on the main stratum.)

\begin{definition}
The {\em cotangent bundle} $T^*\mathcal{M}$ is the vector bundle (existing
by Swan's theorem) over $M$ whose space of sections
is $\Lambda^1(\cM)$.
\end{definition}

\begin{remark}
Note that the elements of $\Lambda^1(\cM)$ are precisely the
forms vanishing on the fibers of the projection $\pi:M\lra\cM$.
\end{remark}

\begin{remark}
Since differential operators on $M$ are just polynomials in vector fields
$X\in\Vect_\cM$ with smooth ($C^\infty(M)$) coefficients,
one can readily show that their interior symbols
are smooth functions on $T^*\mathcal{M}$.
\end{remark}

The cotangent bundles of the strata $\cM_j$ are defined similarly.

\paragraph{The space $C^\infty(\cM)$.}

By definition, the elements of this space are smooth functions on the
blow-up $M$ of $\cM$ that depend
only on $x$ for small $r$ in the
coordinates $(x,r,\om)$ in a neighborhood of each stratum of
nonmaximal dimension.

\subsection{Pseudodifferential operators and symbols}

Here we describe an algebra of zero-order {$\psi$DO} on a stratified
manifold. More precisely, we define families of {$\psi$DO} with
parameter $v$ ranging in a finite-dimensional vector space $V$; by
successive trivial generalizations at each inductive
step, we also obtain {$\psi$DO} smoothly depending on
some additional parameters $x$ and then
{$\psi$DO} parametrized by points of a finite-dimensional vector bundle
over a smooth manifold. Details of these generalizations are left
to the reader.

\paragraph{Negligible families.}

We introduce the space of operator families modulo
which {$\psi$DO} will be defined below. Let $\cM$ be a stratified
manifold (possibly, noncompact). Let $J_\infty(V,\cM)\equiv
J_\infty(\cM)\equiv J_\infty$ be the space of smooth operator
families
\begin{equation}\label{op-sem}
D(v):L^2(\cM)\lra L^2(\cM),\quad v\in V,
\end{equation}
such that all
operators $D(v)$ are compact in $L^2(\cM)$ and satisfy the
estimates
\begin{equation}\label{J0}
\norm{\pd{^\be D(v)}{v^\be}}\le
C_{\be N}(1+\abs{v})^{-N},\qquad\abs{\be},N=0,1,2,\dotsc\,,
\end{equation}
and these conditions hold if we replace $D(v)$ by the product
$$
V_1\dotsm V_pD(v)V_{p+1}\dotsm V_{p+q}
$$
of arbitrary length $p+q\ge 0$. Here $V_1,\dotsc,V_{p+q}$ are smooth vector
fields on $M$ that have the form $V=(V(x),0,\wt
V(x,\om))$ in local coordinates $(x,r,\om)\in\RR^m\times\ov\RR_+\times\Om$
near each stratum, where $V(x)$ is a smooth vector field on the
stratum and $\wt V(x,\om)$ is a vector field, smoothly depending on $x\in X$, with the same
properties on the manifold $\Om$ of lower dimension.

\paragraph{{$\psi$DO} on smooth manifolds.}

We are now in a position to define {$\psi$DO}. Our definition is by induction. We start from the class of families of
{$\psi$DO} on smooth manifolds.

\begin{definition}\label{basis}
An operator family
$$
D(v): L^2(\mathcal{M})\longrightarrow L^2(\mathcal{M})
$$
on a smooth manifold $\cM$ is called a \textit{pseudodifferential
operator with parameter} $v\in V$ if it is a zero-order {$\psi$DO} on
$\cM$ with parameter $v\in V$ in the sense of Agranovich--Vishik.
(These are often called parameter-dependent $\psi$DO.)

The \emph{symbol} (corresponding to the unique stratum of $\cM$) is by
definition the symbol $\si(D)(x,\xi,v)$ in the sense of Agranovich--Vishik.
It is defined on the total space of the vector bundle
$T^*\cM\times V$ over $\cM$ minus the zero section.
\end{definition}

The space of pseudodifferential operators with parameter $v\in V$ on $\cM$
is denoted by $\Psi(V,\cM)$. (In what follows, $V$ will
be suppressed if it is clear from the context or trivial.)

\paragraph{{$\psi$DO} on stratified manifolds.}

Here we define {$\psi$DO} by induction over the length $k$ of the stratified
manifold. Note that we simultaneously define {$\psi$DO} and their symbols.
Definition~\ref{basis} will serve as the inductive statement for $k=0$.

\begin{definition}\label{pdo-induction}
Let $\cM$ be a stratified manifold of length $k>0$. A smooth family of linear
operators
$$
D(v): L^2(\mathcal{M})\longrightarrow L^2(\mathcal{M})
$$
is called a \textit{pseudodifferential operator on} $\cM$ (\textit{with
parameter $v\in V$ in the sense of Agranovich--Vishik}) if the following
conditions hold.
\begin{enumerate}
\item Given $\ph,\psi\in C^\infty(\cM)$ such that
$\supp\ph\cap\supp\psi=\varnothing$, we have $\psi A\ph\in J_\infty$.
\item $D$ is a {$\psi$DO} with parameter on\footnote{The stratum
$\cM_0$ has measure zero in $\cM$. Hence $D$ is automatically interpreted
as an operator on $L^2(\cM\setminus\cM_0)$.} $\cM\setminus\cM_0$; in a
neighborhood $U$ of $\cM_0$, the operator $D$ is representable in the form
\begin{equation}\label{vid}
D=P\BL(x,r,rv,-ir\pd{}x,ir\pd{}r+i\frac{n+1}2\BR),\qquad n=\dim\Omega
\end{equation}
modulo the ideal $J_\infty(V,\cM)$, where $P(x,r,v,\eta,p)\in\Psi(V\times
T^*_x\cM_0\times\RR,\Om)$ is a {$\psi$DO} with parameters on $\Om$,
depending smoothly on additional parameters $x\in\cM_0$ and
$r\in\overline{\RR}_+$; moreover, $P=0$ for $r>r_0$, where $r_0$ is
small.
\end{enumerate}

By definition, the symbols $\si_j(D)$ of $D$ corresponding to the
strata $\cM_j\setminus\cM_{j-1}$, $j>0$, of $\cM$, are its symbols as
an element of
$\Psi(V,\cM\setminus\cM_0)$. The \textit{symbol of} $D$
\textit{corresponding to the stratum} $\cM_0$ is the operator family
\begin{equation}\label{symbol0}
\si_0(D)=P\BL(x,0,rv,r\xi,ir\pd{}r+i\frac{n+1}2\BR):L^2(K_\Om)\lra
L^2(K_\Om)
\end{equation}
parametrized by the points of the bundle $V\times T^*\cM_0$ over $\cM_0$ minus
the zero section.

\textit{The symbols $\si_j(\si_0(D))$ of the symbol} $\si_0(D)$ are just the
symbols $\si_j(P(x,0,v,\eta,p))$ of the {$\psi$DO} $P(x,0,v,\eta,p)$ with
parameters on $\Om$, $j=1,\ldots,k$.
\end{definition}

Let us give some explanations concerning formulas~\eqref{vid}
and~\eqref{symbol0}. One can prove by induction over the length of $\cM$
that a pseudodifferential operator $A\in\Psi(V,\cM)$
satisfies the estimates
\begin{equation}\label{esti-psido}
\norm{\pd{^\a A(v)}{v^\a}}\le
C_\a(1+\abs{v})^{-|\a|},\qquad\abs{\a}=0,1,2,\dotsc\,,
\end{equation}
where all derivatives starting from the first
are compact-valued. Indeed, using this fact for $\cM=\Om$, we see that the
operator family
\begin{equation*}
F(x,t,v,\xi,p)=P(x,e^{-t},ve^{-t},\xi e^{-t},p)
\end{equation*}
satisfies the estimates
\begin{multline}\label{Opa}
\norm{\pd{^{} F(x,t,v,\xi,p)}{x^\a\pa t^l\pa v^\be\pa\xi^\ga\pa p^k}}\le
C_{\a l \be \ga k}(1+\abs{v}+\abs{\xi})^{-|\be|-|\ga|}(1+\abs{p})^{-k},\\
\abs{\a}+l+|\be|+|\ga|+k=0,1,2,\dotsc,
\end{multline}
and the operator family
\begin{equation*}
\wt F(x,t,v,\xi,p)=P(x,0,ve^{-t},\xi e^{-t},p)
\end{equation*}
satisfies the estimates
\begin{multline}\label{Opb}
\norm{\pd{^{}\wt F(x,t,v,\xi,p)}{x^\a\pa t^l\pa v^\be\pa\xi^\ga\pa p^k}}\le
C_{\a l \be \ga k}(e^t+\abs{v}+\abs{\xi})^{-|\be|-|\ga|}(1+\abs{p})^{-k}\\\le
C_{\a l \be \ga
k}(\abs{v}+\abs{\xi})^{-|\be|-|\ga|}(1+\abs{p})^{-k},\quad
\abs{\a}+l+|\be|+|\ga|+k=0,1,2,\dotsc.
\end{multline}
Moreover, both families have compact variation in the parameters
$(v,\xi,p)$. It is not difficult to show that the
operators on the right-hand sides in~\eqref{vid} and~\eqref{symbol0} are
well defined as {$\psi$DO} in the sense of Luke~\cite{Luk1}. Indeed, the
change of variable $r=e^{-t}$ transforms the cone $K_\Om$
into the cylinder $\Om\times\RR$ and the
operator $ir\pa/\pa r$ into $-i\pa/\pa t$. It remains to note that the
operator
$-i\pa/\pa t+i(n+1)/2$ is self-adjoint in the $L^2$ space with the weight $e^{-(n+1)t}$ on the cylinder. This
space is just the image of $L^2(K_\Om)$ under our change of variables.
Therefore, the substitution of this operator as an operator argument
into \eqref{symbol0} is well defined. In
addition, the constructed {$\psi$DO} satisfies the
estimates~\eqref{esti-psido}. This ends the inductive step.

\begin{remark}
Since the cone is noncompact (in $r$), the operator-valued
symbol~\eqref{symbol0} only has \textit{almost compact} variation in
$(\xi,v)$ in the general case. (That is, the
variation becomes compact if we multiply it by a cut-off function
finite in $r$). However, as we shall see shortly, the fiber variation
of the symbol in~\eqref{symbol0} is compact provided
that all symbols $\si_j(\si_0(D))$ are zero,
$j=1,\dotsc,k$.
\end{remark}

\begin{definition}\label{AA}
The \emph{conormal symbol} $\sigma_c(\sigma_0(D))\in
\Psi(\mathcal{M}^\circ_0\times V\times \mathbb{R},\Omega)$ of the family
\eqref{symbol0} is the family $\sigma_c(\sigma_0(D))=P(x,0,v,0,p).$
\end{definition}

\paragraph{Compatibility conditions.}

We have defined the notion of {$\psi$DO} with parameters on a stratified
manifold $\cM$. Such {$\psi$DO} have symbols defined a priori on the
cotangent bundles of the \emph{open strata}
times the parameter space minus the zero section. However, if we
compare the representation~\eqref{vid} with the representation valid in
$U\setminus\cM_0$ by the inductive statement, then we can see that actually
the symbols extend continuously (and smoothly) up to the boundary of the
cotangent bundle and hence are defined on the cotangent
bundles of the corresponding \emph{closed strata}. Moreover, at the points
where a stratum $\cM_j$ meets a stratum $\cM_i$, $j>i$, the following
compatibility conditions hold:
\begin{equation}\label{soglas}
\si_l(\si_j(D))|_{\cM_i}=\si_l(\si_i(D)),\qquad l=j,\dotsc,k.
\end{equation}
Here we write $\si_j(\si_j(D))=\si_j(D)$ for brevity.

\paragraph{Main properties of the calculus of {$\psi$DO}.}

Let us state the main properties of the calculus of
{$\psi$DO} on a compact stratified manifold. The proofs (except for
the proof of Proposition~\ref{dvatri}) are omitted; they can be found
in~\cite{NaSaSt2}.

The set of all symbols~\eqref{symbol0} on $X=\cM_0$ is denoted by
$\Psi({T^*X\times V},K_\Om)$.

\begin{proposition}
$\Psi({T^*X\times V},K_\Om)$ is a local $C^*$-algebra.
\end{proposition}
The norm is given by the supremum of the operator norm over all parameter
values.

\begin{theorem}[Main properties of {$\psi$DO}]\label{svoistva-pdo}
Pseudodifferential operators enjoy the following properties.

\rom{(1)} $\Psi(V,\cM)$ is an algebra with the usual composition of
operators and is a local $C^*$-algebra with respect to the supremum of the
operator norm over the parameter space. Pseudodifferential operators
compactly commute with the operators of multiplication by continuous
functions on $\cM$.

\rom{(2)} The symbol map
\begin{equation}\label{symotobr}
\begin{aligned}
\si: \Psi(V,\cM) &\lra\bigoplus_{j=0}^k \Psi({T^*\cM_j\times
V},K_{\Om_j}),\\
D&\longmapsto (\si_0(D),\dotsc,\si_k(D))
\end{aligned}
\end{equation}
is a local $C^*$-algebra homomorphism and induces an isomorphism
\begin{equation*}
\si:\Psi(V,\cM)\bigm\slash J(V,\cM)\lra\Si(V,\cM)
\subset \bigoplus_{j=0}^k \Psi({T^*\cM_j\times
V},K_{\Om_j})
\end{equation*}
onto the local $C^*$-algebra of symbols subject to the compatibility
conditions~\eqref{soglas}. Here $J(V,\cM)\subset\Psi(V,\cM)$ denotes the
ideal of compact-valued operator families vanishing at infinity.

\rom{(3)} $\Psi(V,\mathcal{M})$ is invariant under diffeomorphisms of $\cM$.
\end{theorem}

\begin{definition}\label{defell}
An operator $D\in\Psi(V,\cM)$ is said to be
\textit{elliptic} if all its symbols $\si_j(D)$, $j=0,\dotsc,k$, are
invertible in the complements of the zero sections of the corresponding
bundles.
\end{definition}
As a corollary of Theorem~\ref{svoistva-pdo}, we obtain the
following assertion
\begin{theorem}\label{ellipticity}
\rom{(1)} On a compact manifold $\cM$, elliptic operators are Fredholm \rom(for
all parameter values\rom).

\rom{(2)} If $V\ne\{0\}$, then an elliptic operator with parameter is
invertible for large $\abs{v}$.
\end{theorem}

The following proposition will be useful.
\begin{proposition}\label{dvatri}
Let $\Si_0\subset \Psi(V,K_\Om)$ be the set of symbols whose symbols
corresponding to the strata of $\Om\times\RR_+$ are zero.
Then each symbol $\si\in\Si_0$ has compact variation in $v$.
\end{proposition}

\begin{proof}
By the definition (see Eq.~\eqref{symbol0}),
$$
\si(v)=P\BL(rv,ir\pd{}r+i\frac{n+1}2\BR):L^2(K_\Om)\lra
L^2(K_\Om)
$$
for an operator function $P(w,p)\in\Psi(V\times \RR,\Om)$.
Theorem~\ref{svoistva-pdo} implies that $P(w,p)\in J(V\times\RR,\Om)$ and
$P(w,p)$ satisfies the estimates~\eqref{esti-psido} with
respect to $(w,p)$.

We have to show that, given $v\ne 0$, $\partial\sigma/\partial v$ is compact.
We have
$$
\frac{\partial \sigma}{\partial v}=r\frac{\partial P}{\partial w}
\BL(rv,ir\pd{}r+i\frac{n+1}2\BR).
$$
The symbol $r\partial P/\partial w(rv,p)$ possesses all estimates
needed for it to define a bounded operator on
$L^2(K_\Omega)$ (cf.\ \eqref{Opb}). Furthermore, it
is compact-valued. It remains to show that it tends to zero as $r\to 0$ or
$r\to\infty$. For $r\to 0$, this is obvious (by virtue of the factor
$r$), and for $r\to \infty$ we use the representation
$$
r\frac{\partial P}{\partial w}(rv,p)=\frac 1{|v|}|rv| \frac{\partial
P}{\partial w}(rv,p)= \frac 1{|v|}\Bigl[|w|\frac{\partial P}{\partial
w}\Bigr]_{w=rv}
$$
and apply the following lemma to $P(w,p)$.
\begin{lemma}
Let $f\in C^2(\RR^n)$ be an operator function such that $f(\xi)\to0$
as $\xi\to\infty$ and $\abs{f''(\xi)}\le C\abs{\xi}^{-2}$ for
large $\xi$. Then $\abs{\xi}\abs{f'(\xi)}\to0$ as $\xi\to\infty$.
\end{lemma}
\end{proof}

\paragraph{Semiclassical quantization.}

Consider the quantization~\eqref{symbol0} with a ``semiclassical''
parameter $h$:
$$
T_h:\Psi({T^*X\times\mathbb{R}},\Omega)\longrightarrow
\Psi({T^*X},K_\Omega)\subset \mathcal{B}(L^2(K_\Omega)).
$$
This map takes a family with parameter on the base of the cone to a family of
operators on the infinite cone and is defined as
\begin{equation}\label{asia}
(T_hD)(\xi):=D\BL(\ovs2r\xi,ih\ovs1{r\pd{}r}+ih(n+1)/2\BR).
\end{equation}
By the same argument as in \cite[Appendix]{NSScS14}, it can be proved that this
quantization is asymptotic in $L^2$, i.e., satisfies the estimates
\begin{equation}\label{7}
T_h(a)T_h(b)=T_h(ab)+o(1),\quad h\to 0,
\end{equation}
in the sense of operator norm.

We shall use this quantization below in the computation of the boundary map in
$K$-theory of {$\psi$DO} algebras.

\section{$\Ell$-groups}
From now on, we fix a compact stratified manifold
$\mathcal{M}$.
\begin{definition}{
Two elliptic operators
$$
D:L^2(\mathcal{M},E)\to L^2(\mathcal{M},F)\text{ and }
D':L^2(\mathcal{M},E')\to L^2(\mathcal{M},F')
$$
acting between sections of vector bundles over the blow-up
$M$ are said to be \emph{stably homotopic} if
there exists a continuous homotopy \footnote{We can and will always
assume that the homotopies are such that all symbols and their
derivatives are continuously homotopic in the corresponding norms.} of elliptic operators
$$
D\oplus 1_{E_0}\sim f^*\bigl(D'\oplus 1_{F_0}\bigr)e^*,
$$
where $E_0,F_0\in\Vect(M)$ are vector bundles and
$$
e:E\oplus E_0\longrightarrow E'\oplus F_0,\qquad f:F'\oplus F_0\longrightarrow F\oplus E_0
$$
are vector bundle isomorphisms. }
\end{definition}

Here ellipticity means the invertibility of symbol components on all
the strata (see Def.~\ref{defell}), and we consider only homotopies of
$\psi$DO preserving ellipticity.

\textbf{Even groups $\Ell_0(\mathcal{M})$}. Stable homotopy is an
equivalence relation on the set of all elliptic pseudodifferential
operators acting between vector bundle sections. Let $\Ell_0(\mathcal{M})$
be the quotient by this equivalence relation. This set is a group with
respect to the direct sum of elliptic
operators. The inverse is defined by the almost inverse operator (that
is, inverse modulo compact operators).

\textbf{Odd groups $\Ell_1(\mathcal{M})$}. We can similarly define the odd
elliptic theory $\Ell_1(\mathcal{M})$ as the group of stable homotopy
classes of elliptic self-adjoint operators. In this case, the stabilization
is done in terms of the operators $\pm Id$.

\begin{remark}
An equivalent definition of the odd $\Ell$-group can be given in terms of
elliptic families on $\mathcal{M}$ parametrized by the circle $\mathbb{S}^1$
modulo constant families.
\end{remark}

\medskip

We consider the following \emph{homotopy
classification problem for elliptic operators:}
compute the group $\Ell_*(\mathcal{M})$.

\medskip

\section{Main theorem}

\paragraph{The map into $K$-homology.}
Let
$$
D:L^2(\mathcal{M},E)\longrightarrow L^2(\mathcal{M},F)
$$
be an elliptic operator. By Theorems~\ref{ellipticity}
and~\ref{svoistva-pdo}, (1), this operator can be
treated as an abstract elliptic operator in the
sense of Atiyah~\cite{Ati4} on $\mathcal{M}$. Thus it
defines an element in the $K$-homology of $\mathcal{M}$. The
corresponding Fredholm module is defined in the standard way. Namely
(cf.~\cite{HiRo1}), if $D$ is self-adjoint (and $E=F$), then consider the
normalization
\begin{equation}\label{oda1}
\mathcal{D}=(P_{\ker D}+D^2)^{-1/2}D:L^2(\mathcal{M},E)\longrightarrow
L^2(\mathcal{M},E),
\end{equation}
where $P_{\ker D}$ is the projection onto the null space
of $D$.

In the general case, we consider the self-adjoint operator
\begin{equation}\label{eva1}
\mathcal{D}=\left(\begin{array}{cc}
0 & \!\!\!\!\!\! D(P_{\ker D}+D^*D)^{-1/2} \\
(P_{\ker D}+D^*D)^{-1/2}D^*\!\!\!\!\!\!& 0
\end{array}
\right):L^2(\mathcal{M},E\oplus F)\rightarrow
L^2(\mathcal{M},E\oplus F),
\end{equation}
which is odd with respect to the
$\mathbb{Z}_2$-grading of $L^2(\mathcal{M},E)\oplus L^2(\mathcal{M},F)$.

\begin{proposition}\label{modul}
1. The operators \eqref{oda1} and \eqref{eva1} define $K$-homology
elements denoted by $[D]\in K_*(\mathcal{M}),$ where $*=1$ in the
self-adjoint case and $*=0$ otherwise.

2. There is a well-defined group homomorphism
$$
\begin{array}{ccc}
\Ell_*(\mathcal{M}) & \stackrel{\varphi}\longrightarrow & K_*(\mathcal{M}),\\
D&\mapsto & [D].
\end{array}
$$

\end{proposition}

\noindent\emph{Proof.} The operators $\mathcal{D}$ in \eqref{oda1} and
\eqref{eva1} are self-adjoint and act in $*$-modules
over the $C^*$-algebra $C(\mathcal{M})$. To end the proof of the first part
of the proposition, it suffices to check that, given $f\in
C(\mathcal{M})$, one has
\begin{equation}\label{kasparov}
[\mathcal{D},f]\in \mathcal{K},\quad (\mathcal{D}^2-1)f\in \mathcal{K},
\end{equation}
where $\mathcal{K}$ is the ideal of compact operators. The compactness
follows easily from the composition formula for pseudodifferential
operators (since $\mathcal{D}$ is a pseudodifferential operator). The
map is well defined, since homotopies of elliptic operators give
continuous homotopies of the corresponding Fredholm modules, i.e., give the
same $K$-homology element. Bundle isomorphisms give degenerate Fredholm
modules. (Recall \cite{HiRo1} that a module is degenerate if all
expressions in \eqref{kasparov} are zero.) \dokaend

\paragraph{The classification theorem.}

The following theorem solves the classification problem on stratified
manifolds.
\begin{theorem}\label{thmain0}
The map
$$
\Ell_*(\mathcal{M}) \stackrel{\varphi}\simeq K_*(\mathcal{M})
$$
that takes each elliptic operator $D$ to the element defined
in Proposition~\emph{\ref{modul}} is an isomorphism.
\end{theorem}
The nondegeneracy of the index pairing $K_0(\mathcal{M})\times
K^0(\mathcal{M})\longrightarrow \mathbb{Z}$ (on the torsion-free
parts of the groups) gives the following assertion.
\begin{corollary}
Two elliptic operators $D_1$ and $D_2$ are stably rationally homotopic if and
only if their indices with coefficients in any vector bundle over $\mathcal{M}$
are equal.
\end{corollary}

We shall obtain Theorem~\ref{thmain0} as a special case of a more general
theorem, which we now state.

\paragraph{The classification of partially elliptic operators.} An operator $D$ on
$\mathcal{M}$ is said to be \emph{elliptic on}
$\mathcal{M}\setminus \mathcal{M}_j$ if the symbol components
$\sigma_k(D),\ldots,\sigma_{j+1}(D)$ are invertible on their domains, i.e.,
everywhere on the complement of the zero section in
$T^*\mathcal{M}_k,\ldots, T^*\mathcal{M}_{j+1}$.

Let $\Ell_*(\mathcal{M}, \mathcal{M}_j)$ be the group of stable homotopy
classes of pseudodifferential operators elliptic on $\mathcal{M}\setminus
\mathcal{M}_j$. We assume that the homotopies are also taken in this
class.

By analogy with the map $\varphi$ in
Proposition~\ref{modul}, we define
$$
\Ell_*(\mathcal{M}, \mathcal{M}_j) \stackrel{\varphi}\longrightarrow
K_*(\mathcal{M}\setminus \mathcal{M}_j).
$$
There is one important difference: to define
$\varphi$ in this case, we replace
$(P_{\ker D}+D^*D)^{-1/2}$ in \eqref{oda1} and \eqref{eva1}
by self-adjoint pseudodifferential operators,
with leading $k-j$ symbol components equal to
$$
(\sigma_k(D)^*\sigma_k(D))^{-1/2},\ldots,
(\sigma_{j+1}(D)^*\sigma_{j+1}(D))^{-1/2}.
$$
Note that the two constructions give the same $K$-homology element for an
operator elliptic on the entire $\mathcal{M}$.

\begin{theorem}\label{thmain}
Given $-1\le j\le k-1$, one has the isomorphism
$$
\Ell_*(\mathcal{M},\mathcal{M}_j) \stackrel{\varphi}\simeq
K_*(\mathcal{M}\setminus \mathcal{M}_j).
$$
\end{theorem}

The proof of Theorem~\ref{thmain} occupies Sections~\ref{pet}--\ref{sem}.
First, in Section~\ref{pet}, $\Ell$-groups are represented as $K$-groups of
certain algebras (noncommutative analog of the Atiyah--Singer difference
construction). This enables us to define exact sequences for $\Ell$-groups.
Then the proof is done inductively over the strata in Sections~\ref{sest}
and \ref{sem}.

\section{Relationship between $\Ell$-groups and $K$-theory}\label{pet}

\textbf{$\Ell$-groups as $K$-groups of $C^*$-algebras}. The embedding
$$
C^\infty(M)\subset \Psi(\mathcal{M})
$$
(the embedding corresponds to the usual action of functions as
multiplication operators) of algebras of scalar operators enables one to
describe pseudodifferential operators in vector bundle sections. Namely, an
arbitrary zero-order {$\psi$DO} acting between spaces of vector bundle
sections can be represented as
$$
D':\im P\longrightarrow \im Q,
$$
where $P$ and $Q$ are matrix projections ($P^2=P,Q^2=Q$) with entries in
$C^\infty(M)$ and $D'$ is a matrix operator with entries
in $\Psi(\mathcal{M})$.

Denote by
$$
\Sigma(\mathcal{M}\setminus \mathcal{M}_j)\stackrel{\rm def}=\im
(\sigma_k,\ldots,\sigma_{j+1})\subset \bigoplus_{l\ge j+1}
C^\infty(S^*M_l,\mathcal{B}(L^2(K_{\Omega_l})))
$$
the algebra of leading $k-j$ components of the principal symbol.

Theorem~4 in~\cite{Sav8} gives isomorphisms
\begin{equation}\label{sa1}
\Ell_*(\mathcal{M}, \mathcal{M}_j)\stackrel\chi\simeq
K_*(\Con(C^\infty(M)\stackrel{f}\to \Sigma(\mathcal{M}\setminus
\mathcal{M}_j)))
\end{equation}
of the $\Ell$-groups and the $K$-groups of special local
$C^*$-algebras. Here
$$
f:C^\infty(M)\longrightarrow \Sigma(\mathcal{M}\setminus \mathcal{M}_j)
$$
is the embedding taking a smooth function on the blow-up $M$ to the
symbol of multiplication operator by this function, and
$$
\Con(A\stackrel{f}\longrightarrow B)=\Bigl\{(a,b(t))\in A\oplus C_0([0,1),B)\;|\;
f(a)=b(0)\Bigr\}
$$
is the mapping cone of the algebra homomorphism $f:A\to B$.

In most cases, one can eliminate the mapping cone
from Eq.~\eqref{sa1} by using the following lemma.
\begin{lemma}
The $K$-group of the mapping cone can be represented as
$$
K_{*+1}(\Con(C^\infty(M)\to \Sigma(\mathcal{M}\setminus \mathcal{M}_j)))\simeq
K_{*}(\Sigma(\mathcal{M}\setminus \mathcal{M}_j))/K_{*}(C^\infty(M)),
$$
provided that $M$ has a nonsingular vector field. \rom(This condition
is satisfied if, say, $M$ has no components with empty
boundary.\rom)
\end{lemma}
\noindent \emph{Proof.} A vector field $M\to S^*M$ defines
a section $\Sigma (\mathcal{M}\setminus \mathcal{M}_j)\to C^\infty(M)$.
Thus the mapping cone exact sequence (induced by the embedding of
algebras) splits. The splitting gives the desired isomorphism. \dokaend

\begin{remark}
In the odd case, the composition of the latter isomorphism with $\chi$
shows that, modulo stable homotopy, elliptic self-adjoint
operators are isomorphic to symbols-projections modulo
projections onto bundle sections (cf.~\cite{APS3}).
\end{remark}

{\bf Exact sequence of the pair in $\Ell$-theory}. Let us
construct the exact sequence corresponding in
elliptic theory to the pair
$$
\mathcal{M}_j\setminus \mathcal{M}_{j-1} \subset \mathcal{M}\setminus \mathcal{M}_{j-1}.
$$
We define the maps in the desired sequence in
$K$-theoretic terms. To this end, consider the commutative
diagram
\begin{equation}
\label{dia1}
\begin{array}{rcccccl}
0\rightarrow & \ker (\sigma_k,\ldots,\sigma_{j+1}) & \longrightarrow &
\Sigma(\mathcal{M}\setminus \mathcal{M}_{j-1}) & \stackrel{(\sigma_k,\ldots,\sigma_{j+1})}
\longrightarrow & \Sigma(\mathcal{M}\setminus \mathcal{M}_{j})
&
\to 0\\
& \uparrow & & \uparrow & & \uparrow \\
0\rightarrow & 0 & \longrightarrow & C^\infty(M) & = & C^\infty(M) & \to
0
\end{array}
\end{equation}
with exact rows. The ideal $\ker (\sigma_k,\ldots,\sigma_{j+1})$ of
symbols with leading $k-j$ components equal to zero is denoted for brevity
by $\Sigma_0$. The diagram induces the exact sequence
\begin{equation}\label{conu1}
0\to S\Sigma_0 \rightarrow \Con(C^\infty(M)\to \Sigma(\mathcal{M}\setminus
\mathcal{M}_{j-1})) \rightarrow \Con(C^\infty(M)\to
\Sigma(\mathcal{M}\setminus \mathcal{M}_{j})) \to 0
\end{equation}
of the mapping cones of the vertical embeddings. Here
$S\Sigma_0=C_0((0,1),\Sigma_0)$ stands for suspension.

The $K$-groups of the mapping cones in \eqref{conu1} classify
elliptic operators on $\mathcal{M}\setminus \mathcal{M}_{j-1}$ and
$\mathcal{M}\setminus \mathcal{M}_{j}$,
respectively (see~Eq.~\eqref{sa1}). It turns out that the $K$-groups
of the ideal also classify elliptic operators.
\begin{lemma}
The $K$-groups of $S\Sigma_0$ classify operators that are elliptic on
$\mathcal{M}\setminus \mathcal{M}_{j-1}$ and have symbols
$(\sigma_k,\ldots,\sigma_{j+1})$ induced by constant functions in
$C^\infty(M)$.
\end{lemma}
\begin{proof}
To be definite, we consider $K_0(S\Sigma_0)$. Here
$K_0(S\Sigma_0)=K_1(\Sigma_0)=K_1(\Sigma_0^+)$, where $\Sigma_0^+$ denotes
the algebra with identity adjoined. Now it follows from
the definition of $K_1$ and $\Sigma_0$ that $K_1(\Sigma_0^+)$ classifies
elliptic operators described in the lemma.
\end{proof}

Using this lemma, we can define a map (cf.\
Proposition~\ref{modul})
\begin{equation}\label{quant1}
\varphi: K_*(S\Sigma_0)\to K_*(\mathcal{M}_j\setminus
\mathcal{M}_{j-1})
\end{equation}
that restricts elliptic operators to a neighborhood $U$
of the subset $\mathcal{M}_j\setminus \mathcal{M}_{j-1}$. We assume that
$U$ has the structure of a bundle $\pi:U\to
\mathcal{M}_j\setminus \mathcal{M}_{j-1}$ with conical fiber. Note
that the structure of a $C_0(\mathcal{M}_j\setminus
\mathcal{M}_{j-1})$-module
$$
C_0(\mathcal{M}_j\setminus \mathcal{M}_{j-1})\stackrel{\pi^*}\longrightarrow
C(U)\to \mathcal{B}(L^2(U))
$$
on the corresponding $L^2$-spaces where the operator acts is obtained by the pull-back from the base of the bundle. The
restriction of the operator is well defined, since in the complement of
an arbitrarily small neighborhood of $\mathcal{M}_j\setminus
\mathcal{M}_{j-1}$ the operator acts as the multiplication by some
function.

We shall show in Lemma~\ref{lem2} that
\begin{equation}\label{excision1}
K_*(S\Sigma_0)\simeq K^{*}_c(T^*M_j)=\Ell_*(\mathcal{M}_{j},
\mathcal{M}_{j-1}).
\end{equation}
Therefore, we can replace all $K$-groups by $\Ell$-groups in the
exact sequence induced in $K$-theory
by~\eqref{conu1}. As a result, we obtain the periodic six-term
exact sequence relating the $\Ell$-groups:
\begin{equation}\label{hexa}
\xymatrix{ & \Ell_0\left(\mathcal{M}_j, \mathcal{M}_{j-1}\right)\ar[r] &
\Ell_0\left(\mathcal{M}, \mathcal{M}_{j-1}\right)
\ar[rd] \\
\Ell_1\left(\mathcal{M}, \mathcal{M}_{j} \right)\ar[ur]^\partial & & &
\Ell_0\left(\mathcal{M}, \mathcal{M}_{j} \right)\ar[dl]^\partial \\
& \Ell_1\left(\mathcal{M}, \mathcal{M}_{j-1}\right)\ar[lu] & \Ell_1\left(\mathcal{M}_j, \mathcal{M}_{j-1}\right).\ar[l] }
\end{equation}
We do not dwell upon the question of defining all maps in
this sequence directly in terms of elliptic operators.

\section{Induction}\label{sest}

For $j$ ranging from $k$ down to $-1$, we shall prove by induction that the map $\varphi$
in Theorem~\ref{thmain} is an isomorphism on the set
$\mathcal{M}\setminus \mathcal{M}_{j}$.

For $j=k$, this is obvious. Let us prove the inductive step:
the isomorphism for $j+1$ implies an isomorphism for
$j$.

The pair $\mathcal{M}_j\setminus \mathcal{M}_{j-1}\subset \mathcal{M}\setminus
\mathcal{M}_{j-1}$ gives the diagram
\begin{equation}\label{inda}
\begin{array}{cccc}
\ldots\to\Ell_{*+1}(\mathcal{M}, \mathcal{M}_{j}) \stackrel\partial\to
&\hspace{-3mm} \Ell_*(\mathcal{M}_j, \mathcal{M}_{j-1}) \to
&\hspace{-3mm}\Ell_*(\mathcal{M}, \mathcal{M}_{j-1}) \to &\hspace{-3mm}
\Ell_*(\mathcal{M}, \mathcal{M}_{j}) \stackrel\partial\to\ldots
\\
\downarrow & \downarrow & \downarrow & \downarrow \\
\ldots\to K_{*+1}(\mathcal{M}\setminus \mathcal{M}_{j}) \stackrel\partial\to &
\hspace{-3mm}K_*(\mathcal{M}_j\setminus \mathcal{M}_{j-1}) \to
&\hspace{-3mm}K_*(\mathcal{M}\setminus \mathcal{M}_{j-1}) \to &
\hspace{-3mm}K_*(\mathcal{M}\setminus \mathcal{M}_{j})
\stackrel\partial\to\ldots
\end{array}
\end{equation}
Note that the vertical maps are defined
in Proposition~\ref{modul} and Eqs.~\eqref{quant1}
and~\eqref{excision1}. Once we prove that
the diagram commutes, the inductive assumption, in conjuntion
with the $5$-lemma, will imply that the
vertical map ranging in $K_*(\mathcal{M}\setminus \mathcal{M}_{j-1})$
is an isomorphism, as desired. Let us establish the
commutativity.

The commutativity of the square
$$
\begin{array}{ccc}
\Ell_*(\mathcal{M}, \mathcal{M}_{j-1}) &\longrightarrow & \Ell_*(\mathcal{M},
\mathcal{M}_{j})
\\
\downarrow & & \downarrow\\
K_*(\mathcal{M}\setminus \mathcal{M}_{j-1}) & \longrightarrow &
K_*(\mathcal{M}\setminus \mathcal{M}_{j})
\end{array}
$$
is obvious, since the horizontal maps are induced by forgetting some of the
structures.

The commutativity of the square
$$
\begin{array}{ccc}
\Ell_*(\mathcal{M}_j, \mathcal{M}_{j-1}) &\longrightarrow &
\Ell_*(\mathcal{M},\mathcal{M}_{j-1})
\\
\downarrow & & \downarrow\\
K_*(\mathcal{M}_j\setminus \mathcal{M}_{j-1}) & \longrightarrow &
K_*(\mathcal{M}\setminus\mathcal{M}_{j-1}),
\end{array}
$$
corresponding to the embedding $i:\mathcal{M}_j\setminus
\mathcal{M}_{j-1}\longrightarrow \mathcal{M}\setminus\mathcal{M}_{j-1}$, is
also easy to prove. Indeed, consider the composition of the maps through
the top right corner of the square: it takes an operator elliptic on
$\mathcal{M}\setminus\mathcal{M}_{j-1}$ and equal to a
multiplication operator in the complement of a neighborhood $U$ of
$\mathcal{M}_j\setminus \mathcal{M}_{j-1}$ to the same operator
with the natural $C_0(\mathcal{M}\setminus\mathcal{M}_{j-1})$-module
structure on the spaces between which the operator
acts. We restrict the operator to a neighborhood of $\mathcal{M}_j\setminus
\mathcal{M}_{j-1}$ (this does not change the $K$-homology element, since
the corresponding Fredholm module changes by a degenerate module)
and make a homotopy of the module structure to the composition
$$
C_0(\mathcal{M}\setminus \mathcal{M}_{j-1})\stackrel{i^*}\to
C_0(\mathcal{M}_j\setminus \mathcal{M}_{j-1})\stackrel{\pi^*}\to C(U)\to
\mathcal{B}(L^2(U)),
$$
where $\pi:U\to \mathcal{M}_j\setminus \mathcal{M}_{j-1}$ is the
projection. This gives the very same element as the composition
of arrows through the bottom left corner of the square. This proves
the commutativity of the square.

To justify the inductive step, it remains to show the commutativity
of the squares containing the boundary maps
in~\eqref{inda}. This is the most technically involved part of the proof.
We do this in the next section.

\section{Comparison of boundary maps}\label{sem}

In this section, we prove the commutativity of the squares
containing the boundary maps in~\eqref{inda}.
The proof is carried out in terms of $K$-groups of
symbol algebras. The proof is rather long; therefore, here we
list the main steps:
\begin{itemize}
\item We show that the $K$-homology boundary map is the
composition of the restriction to the boundary of the
stratum and the direct image (see Subsec.~1 below).

\item The $K$-theory boundary map is also the composition of the
restriction map and a certain boundary map $\partial'$ related to the
algebra $\Psi(T^*M_j,K_\Omega)$ of symbols on $M_j$ (Subsec.~2).
Hence the comparison of the boundary maps in $K$-theory and
$K$-homology is reduced to $M_j$.

\item Unfortunately, the boundary map $\partial'$ is not so easy to work
with. Therefore, we replace $\Psi(T^*M_j,K_\Omega)$ by
the simpler algebra $\Psi(T^*M_j\times\mathbb{R},\Omega)$ of families
with parameters. It turns out that the boundary map $\partial''$ for the
latter algebra is easy do compute. More precisely, in Subsec.~3 we
define an asymptotic homomorphism of one algebra into the other, and
in Subsec.~4 we show that the asymptotic homomorphism
induces an isomorphism on $K$-groups.

\item The boundary map $\partial''$ is computed in Subsec.~5 in terms of
the families index of an elliptic family with parameters. Then the
compatibility of the boundary map $\partial''$ and the direct image in
$K$-homology is verified in Theorem~\ref{luke2} in
Subsec.~5. This finishes the proof of the
commutativity in~\eqref{inda}.
\end{itemize}

Let us now give the details of the proof.

\paragraph{1. The boundary map in $K$-homology in \eqref{inda}.}
We start with some notation. Let $U$ be a
neighborhood of $\mathcal{M}_j\setminus
\mathcal{M}_{j-1}=\mathcal{M}_j^\circ$
fibered with conical fiber over the stratum, and
let $\pi: U\to \mathcal{M}_j\setminus \mathcal{M}_{j-1}$ be the
corresponding projection. The bundle of cone bases is
denoted by
$$
\pi':\Xi\longrightarrow \mathcal{M}_j\setminus \mathcal{M}_{j-1}.
$$
The typical fiber is $\Omega$. One obviously has $U\setminus
\mathcal{M}_j^\circ\simeq \mathbb{R}_+\times \Xi$.

We claim that the boundary map $\partial:K_*(\mathcal{M}\setminus
\mathcal{M}_j)\to K_{*+1}(\mathcal{M}_j\setminus \mathcal{M}_{j-1})$ is equal
to the composition
\begin{equation}\label{part1}
K_*(\mathcal{M}\setminus \mathcal{M}_{j})\longrightarrow K_*(U\setminus
\mathcal{M}_j^\circ)= K_{*+1}(\Xi) \stackrel{\pi'_*}\longrightarrow
K_{*+1}(\mathcal{M}_j\setminus \mathcal{M}_{j-1})
\end{equation}
of the restriction to $U\setminus \mathcal{M}_j^\circ$, the
periodicity isomorphism, and the direct image map. This
decomposition is easy to obtain by using the fact that the boundary map is natural.

\paragraph{2. Reduction to the boundary.} Let us start the computation of the
boundary map in $K$-theory of algebras. First,
being natural, the boundary map
$$
\partial:K_*(\Con(C^\infty(M)\to \Sigma(\mathcal{M}\setminus \mathcal{M}_{j})))\to
K_{*+1}(S\Sigma_0)
$$
is equal to the composition
\begin{multline}\label{part2a}
K_*(\Con(C^\infty(M)\to \Sigma(\mathcal{M}\setminus \mathcal{M}_{j})))\longrightarrow K_*(\Con(C^\infty(\partial_j M)\to
\Sigma_{M_j}))\stackrel{\partial'}\longrightarrow K_{*+1}(S\Sigma_0)
\end{multline}
of the restriction of symbols to $M_j$ and the boundary map
$\partial'$, where $\Sigma_{M_j}\equiv\Sigma(T^*M_j,\Omega)$ denotes the
algebra of symbols of $\psi$DO on $\Omega$ with parameters in
$T^*M_j$ and $\partial_j M\subset\partial M$ is the
closure of $\pi^{-1}\mathcal{M}^\circ_j$.\footnote{One can readily see that $\partial_j M$ is a manifold fibered over
$M_j$ with fiber isomorphic to the blow-up of $\Omega$.}

\paragraph{3. The asymptotic homomorphism.} Recall that
$\Psi({T^*M_j},K_\Omega)$ is the algebra of $j$th symbols on $M_j$. To
compute the boundary map $\partial'$ in \eqref{part2a}, we replace
$\Psi({T^*M_j},K_\Omega)$ (preserving its $K$-group) by an algebra of
pseudodifferential operators with parameters,
for which the boundary map is simpler.

To this end, consider the map (see~\eqref{asia})
$$
T_h:\Psi({T^*M_j\times \mathbb{R}},\Omega) \longrightarrow
\Psi({T^*M_j},K_\Omega), \quad h\in (0,1],
$$
\begin{equation}\label{quasi}
(T_h u)(\xi)=u\left(\stackrel{2} r\xi,ih \stackrel 1
{r\frac\partial{\partial r}}+ih(n+1)/2\right), \quad (\xi,p)\in
T^*M_j\times \mathbb{R},
\end{equation}
where $\Psi({T^*M_j\times \mathbb{R}},\Omega) $ is the algebra of smooth
families of {$\psi$DO} on the fibers $\Omega$ with parameters in
$T^*M_j\times \mathbb{R}$. As $ h\to 0$, we have
\begin{equation}\label{asia2}
T_h (ab)= T_h(a)T_h(b)+o(1),\;\;\; \left(T_h(a)\right)^*=T_h(a^*)+o(1),
\end{equation}
where $a,b\in \Psi({T^*M_j\times \mathbb{R}},\Omega)$ are arbitrary
and $o(1)$ is understood in the uniform norm.

Note that this semiclassical quantization is a special case of the
so-called asymptotic homomorphisms, playing an important role in
$C^*$-algebra theory~\cite{CoHi4,Hig2,Man1}. In particular,
Eqs.~\eqref{asia2} imply that $T_h$ induces a
$K$-group homomorphism
$$
T:K_*({\Psi({T^*M_j\times \mathbb{R}},\Omega)}) \to
K_*(\Psi({T^*M_j},K_\Omega)).
$$

Consider the commutative diagram
\begin{equation}\label{dia4}
\begin{array}{rcccccl}
0\to& J({T^*M_j\times \mathbb{R}},\Omega) & \to & \Psi({T^*M_j\times
\mathbb{R}},\Omega) & \to
& \Sigma_{M_j} & \to 0\\
& \;\;\;\downarrow T_h & & \downarrow T_h & & \downarrow t_h\\
0\to &\Sigma_0 & \to & \Psi({T^*M_j},K_\Omega) & \to & \Sigma_{M_j} &
\to 0,
\end{array}
\end{equation}
where $t_h$ is the induced map on the symbols. (Actually, it is
an algebra isomorphism.)

The algebra $C^\infty(\partial_j M)$ is a subalgebra in
$\Psi({T^*M_j\times \mathbb{R}},\Omega),\Sigma_{M_j}$ and
$\Psi({T^*M_j},K_\Omega)$. The diagram of the mapping cones of these
embeddings gives the square
\begin{equation}\label{quadrat}
\begin{array}{ccc}
K_*(\Con(C^\infty(\partial_j M)\to
\Sigma_{M_j})) & \stackrel{\partial''}\longrightarrow & K_{*}(C_0(T^*M_j\times \mathbb{R}))\\
\parallel & &\;\;\;\downarrow T \\
K_*(\Con(C^\infty(\partial_j M)\to \Sigma_{M_j})) &
\stackrel{\partial'}\longrightarrow& K_{*+1}(S\Sigma_0).
\end{array}
\end{equation}
(The horizontal maps are just the boundary maps in the
corresponding sequences.) Here we have used the fact that
$$
K_*(J({T^*M_j\times \mathbb{R}},\Omega))\simeq K_*(C_0(T^*M_j\times
\mathbb{R})).
$$
This isomorphism is induced by the embedding $J({T^*M_j\times
\mathbb{R}},\Omega)\subset C_0(T^*M_j\times \mathbb{R},\mathcal{K})$ of
a local $C^*$-algebra into its closure.

The commutativity in~\eqref{quadrat} follows, since the $K$-theory boundary map is natural with respect
to asymptotic homomorphisms (e.g., see~\cite{Con1}).

\paragraph{4. The map $T:K_*(C_0(T^*M_j\times \mathbb{R}))\longrightarrow
K_*(\Sigma_0)$ is an isomorphism.} First, we compute $K_{*}(S\Sigma_0)$. To
this end, consider the short exact sequence
\begin{equation}
\label{solv2} 0\to \ker \sigma_c\longrightarrow \Sigma_0
\stackrel{\sigma_c}\longrightarrow J(M_j\times \mathbb{R},\Omega) \to 0
\end{equation}
of local $C^*$-algebras. (Here $\sigma_c$ is the conormal
symbol map, see Def.~\ref{AA}.) The kernel
$\ker\sigma_c$ is formed by families of compact
operators.\footnote{Just as in the theory of
operators on manifolds with isolated conical singularities, a family in
$\Sigma_0$ is compact if and only if its conormal symbol is
zero.} Since the embeddings
$$
\ker\sigma_c\subset C(S^*M_j,\cK(K_\Omega))\quad \text{and} \quad
J(M_j\times\mathbb{R},\Omega)\subset C_0(M_j\times
\mathbb{R},\mathcal{K}(\Omega))
$$
induce isomorphisms in $K$-theory, the $K$-theory long exact sequence
induced by \eqref{solv2} can be written as the top row of the diagram
\begin{equation}
\label{dia2}
\begin{array}{ccccccc}
\rightarrow & K^*(S^*M_j) & \rightarrow & K_*(\Sigma_0) &
\rightarrow & K^*_c(M_j\times \mathbb{R}) & \rightarrow\\
& \parallel & & \downarrow {L} & & \parallel \\
\rightarrow & K^*(S^*M_j) & \rightarrow & K^{*+1}_c(T^*M_j)&
\rightarrow & K^*_c(M_j\times \mathbb{R}) & \rightarrow,
\end{array}
\end{equation}
where the bottom row is the sequence of topological $K$-groups of the pair
$S^*M_j\subset B^*M_j$ of the unit sphere and ball bundles in
$T^*M_j$ and the map $L$ is the difference construction
for {$\psi$DO} with operator-valued symbols in the sense of Luke (see
\cite{Luk1,NSScS15}). Recall that this map is defined as follows: given
$\sigma_j\in\Sigma_0^+$, we treat it as an
operator-valued function on $S^*M_j$ of compact variation in the fibers of
$S^*_xM_j$ (see Proposition~\ref{dvatri}); if this function is invertible,
then $L$ takes $[\sigma_j]\in K_1(\Sigma_0)$ to the index
\begin{equation}\label{assa}
L[\sigma_j]:=\ind \widetilde{\sigma_j} \in K_c(T^*M_j)
\end{equation}
of an extension $\widetilde{\sigma_j}$ of $\sigma_j$ to the unit ball
bundle in $T^*M_j$ preserving the compact variation property in the
fibers. (The extension is Fredholm in $B^*M_j$ and invertible
on $S^*M_j$; thus its index is an element of the above
mentioned $K$-group with compact supports.) The
index~\eqref{assa} is independent of the choice of the extension. For
even $K$-groups, the map is defined in a similar way.

\begin{lemma}\label{lem2}
The diagram~\eqref{dia2} commutes, and hence
$L$ is an isomorphism.
\end{lemma}

\begin{proof}
1. The commutativity in
$$
\begin{array}{ccc}
K^*(S^*M_j) & \longrightarrow & K_*(\Sigma_0) \\
\parallel & & \downarrow {L}\\
K^*(S^*M_j) & \longrightarrow & K^{*+1}(T^*M_j)
\end{array}
$$
follows from the fact that for finite-dimensional symbols $L$ coincides
with the Atiyah--Singer difference construction.

2. Consider the square
$$
\begin{array}{ccc}
K_*(\Sigma_0) & \stackrel{\sigma_c}\longrightarrow & K^*_c(M_j\times \mathbb{R}) \\
\downarrow {L} & & \parallel\\
K^{*+1}_c(T^*M_j)& \stackrel{j^*}\longrightarrow & K^{*+1}(M_j),
\end{array}
$$
where $j:M_j\to T^*M_j$ is the embedding of the zero section. Its
commutativity follows from the index formula (e.g.,
see \cite{NSScS3})
\begin{equation}\label{indf1}
\beta \ind D_y=\ind \sigma_c(D_y)\in K^1_c(Y\times \mathbb{R})
\end{equation}
for an elliptic family $D_y,$ $y\in Y$ with unit interior symbol on the
infinite cone. Here $Y$ is a compact parameter space, and
$\beta$ is the periodicity isomorphism $K(Y)\simeq K^1_c(Y\times
\mathbb{R})$.

Indeed, given $a\in K_1(\Sigma_0)$, the elements $j^*L(a)$ and
$\sigma_c(a)$ are, respectively, the
left- and the right- hand side in~\eqref{indf1}. The
case of elements in the $K_0$-group can be verified if
we first consider the suspension and then apply~\eqref{indf1}.

3. Now consider the square
$$
\begin{array}{ccc}
K^{*+1}_c(M_j\times \mathbb{R}) & \stackrel\partial \longrightarrow & K^*(S^*M_j) \\
\parallel & & \parallel\\
K^*(M_j)& \stackrel{p^*}\longrightarrow & K^*(S^*M_j)
\end{array},
$$
where $p:S^*M_j\to M_j$ is the natural projection. Its commutativity
follows from the index formula. More precisely, let us represent $a\in
K^1_c(M_j\times \mathbb{R})$ by an invertible family of conormal symbols
with unit principal symbol. On the one hand, going through the left bottom
corner of the square, we obtain $p^*\ind a\in K^0(S^*M_j)$.
On the other hand, going through the upper right corner, we obtain
$$
\partial a=p^*\ind \widehat{a},
$$
where $\widehat{a}$ denotes the operator family on $K_\Omega$ with unit
interior symbol and with conormal symbol $a$. Note that the latter
equality follows from the fact that $\partial$ takes an invertible symbol
to the index of the corresponding operator. Applying~\eqref{indf1}, we
obtain the desired relation $p^*\ind a=p^*\ind
\widehat{a}\in K^0(S^*M_j)$. For the case of elements in $K_0$ one should
first use suspension.
\end{proof}

\begin{lemma}\label{semdva}
The map
$$
T: K_*(C_0(T^*M_j\times \mathbb{R}))\longrightarrow K_*(\Sigma_0)
$$
is the inverse of the
isomorphism $L$ in~\eqref{dia2}.
\end{lemma}
\noindent \emph{Proof.} To be definite, we consider the map
$$T: K_1(C_0(T^*M_j\times \mathbb{R}))\longrightarrow
K_1(\Sigma_0).$$

Let us prove that $LT$ is the identity map.

1. Given a symbol $u(\xi,p)\in (1+J(T^*M_j\times\RR,\Omega))$ invertible
for $(\xi,p)\in T^*M_j\times \mathbb{R}$ and
equal to the identity in the complement of a compact set,
we obtain
$$
LT [u]=\ind u\left(\stackrel{2} r\xi,ih
\stackrel 1 {r\frac\partial{\partial r}}+ih(n+1)/2\right)\in
K_c^0(T^*M_j)
$$
(for $h$ small), where $[u]\in K_c^1(T^*M_j\times \mathbb{R})$. The index
is well defined, since $(T_h u)(\xi)$ has compact
variation on the fibers in $T^*M_j\setminus \mathbf{0}$ and
is invertible whenever $\xi\ne 0$ (see Proposition~\ref{dvatri} and
formula \eqref{7}, respectively).

2. Let $\overline{T_h u}$ be the family of Fredholm
operators parametrized by $T^*M_j\setminus \mathbf{0}$ and equal to $T_h u$
for $|\xi|<1$ and to
$$
u\left((\stackrel 2 r+|\xi|-1)\xi,ih \stackrel 1 {r\frac\partial{\partial
r}}+ih(n+1)/2\right)
$$
otherwise.

For small $h$, this family is
invertible for all $\xi$. (This follows from the boundedness of the
support of $1-u$ and since~\eqref{7} is uniform
for $\const>\lambda\ge 0$ if in~\eqref{asia} we
replace $r\xi$ by $(r+\lambda)\xi$). We obtain, by
construction,
$$
\ind T_h u=\ind \overline{T_h u}.
$$
However, $\overline{T_h u}$ is the family of
identity operators for large $\xi$. Hence its index
can be calculated by~\eqref{indf1} and is equal to the index of the family
of conormal symbols $u(\xi,p)$ (modulo the Bott periodicity
isomorphism); i.e., it actually gives the original element
$$
LT [u]=[u].
$$
\dokaend

\paragraph{5. Comparison of boundary maps.} Let us compare the obtained
expressions for the boundary maps in $K$-homology and $K$-theory (see
\eqref{part1} and \eqref{part2a}, \eqref{quadrat}). We omit the restriction
maps from Subsecs.~1 and~2. Consider the diagram (here and below,
$I=(0,1)$)
\begin{equation}\label{dia5}
\begin{array}{ccc}
K_{*}(\Xi\times I) & \stackrel{\pi'_*}\longrightarrow &
K_{*}(M^\circ_j\times I)\\
\varphi\uparrow & & \uparrow \varphi\\
K_*(\Con(C^\infty(\partial_j M)\to
\Sigma_{M_j})) & \stackrel{\partial''}\longrightarrow & K^{*}_c(T^*M_j\times \mathbb{R})\\
\parallel & &\;\;\;\uparrow L \\
K_*(\Con(C^\infty(\partial_j M)\to \Sigma_{M_j})) &
\stackrel{\partial'}\longrightarrow& K_{*}(\Sigma_0),
\end{array}
\end{equation}
containing the boundary maps in the middle and bottom rows. All
maps in the diagram have already been defined except for $\varphi$ in
the left column. To define it, we note that $K_*(\Con(C^\infty(\partial_j
M)\to \Sigma_{M_j}))$ classifies elliptic families $\sigma(x,\xi,p)$ on
$\Omega$ with parameters $(x,\xi,p)\in T^*M_j\times \mathbb{R}$. Such a
family defines an operator on the product $\Xi\times I$ by the formula
$$
\sigma\left(x,-i\frac\partial{\partial x},-i\frac\partial{\partial t}\right).
$$
This operator gives an element in $K_{*}(\Xi\times I)$ provided that
$\sigma(x,\xi,p)$ is elliptic.

The bottom square in~\eqref{dia5} is isomorphic to \eqref{quadrat} by
Lemma~\ref{semdva}. Hence it commutes. We claim that the top square
is also commutative. Indeed, consider any element $z\in K_*(\Con(C^\infty(\partial_j
M)\to \Sigma_{M_j}))$ defined by an elliptic family $\sigma(x,\xi,p)$.
Ellipticity implies the Fredholm property, and
$\partial''$ takes $z$ simply to the index of the family with parameters in
$T^*M_j\times \mathbb{R}$. On the other hand, $\pi'_* z$ corresponds to
the elliptic operator
$$
\sigma\left(x,-i\frac\partial{\partial x},-i\frac\partial{\partial t}\right)
$$
in $C_0(M^\circ_j\times I)$-modules. The last two elements
actually coincide. This is a consequence of the
following general result.

\begin{theorem}\label{luke2}
Let $p(x,\xi)$ be an operator-valued symbol elliptic in the sense of Luke
\emph{\cite{Luk1}} on a compact manifold $X$ with corners. Then

\begin{equation}\label{luuk}
\left[p\left(x,-i\frac\partial{\partial x}\right)\right]=\varphi\left(\ind
p(x,\xi)\right)\in K_*(X^\circ),
\end{equation}
where $X^\circ$ is the interior of the manifold, square brackets denote an
element in $K$-homo\-logy, and $\varphi:K_c^*(T^*X)\to K_*(X^\circ)$ is the
Poincar\'e isomorphism on manifolds with corners \rom(e.g.,
see~\emph{\cite{MePi2}}\rom).
\end{theorem}
The proof of this theorem is given in the supplement.

Thus Theorem \ref{luke2} proves that \eqref{dia5} commutes.
Hence we also have the commutativity of the squares in \eqref{inda}
containing the boundary maps.

This concludes the proof of Theorem~\ref{thmain}.

\section{Applications}\label{vosem}

\textbf{A topological obstruction to Fredholm property}. Let
$\mathcal{M}\supset X$ be a stratified pair. It is of
interest to find conditions under which an operator
elliptic in $\mathcal{M}\setminus X$ can be
transformed into an operator elliptic in
$\mathcal{M}$ without changing the components of the principal symbol on
the set $\mathcal{M}\setminus X$. This question is similar to the
Atiyah--Bott problem of determining the topological conditions on the
symbol on a smooth manifold with boundary
under which there exists a Fredholm boundary condition for the
corresponding operator.

We shall answer a similar question for the elements of $\Ell$-groups in the
case of an arbitrary stratification. To this end, consider the diagram
$$
\xymatrix{
\Ell(\mathcal{M}) \ar[r] \ar[d]^\varphi_{\simeq}& \Ell(\mathcal{M}, X)
\ar[d]_\varphi^\simeq\\
K_0(\mathcal{M}) \ar[r] & K_0(\mathcal{M}\setminus X) \ar[r]^{\partial} &
K_1(X).
}
$$
It obviously commutes, since we deal with forgetful maps. Therefore, the
nonvanishing of $\partial \varphi(a)$ is a necessary and sufficient
condition for the existence of a lifting of $a\in \Ell(\mathcal{M}, X)$ to
$\Ell(\mathcal{M})$.

The boundary map in $K$-homology plays a similar role in other problems
(see Baum--Douglas~\cite{BaDo4}, Roe~\cite{Roe2}, and
Monthubert~\cite{Mon2}).

Note that the equation $\partial \varphi(a)=0$ is a
condition on the interior symbol alone
(a finite-dimensional condition) if $X=\mathcal{M}_{k-1}$ is the set
of all singularities of $\mathcal{M}=\mathcal{M}_k$.\vspace{2mm}

\noindent\textbf{Cobordism invariance of the index}. Let us give a
generalization of the usual cobordism invariance of the index of Dirac
operators. Suppose that $X$ is a smooth stratum. Then we have
the commutative diagram
$$
\xymatrix{
\Ell_1(\mathcal{M}, X) \ar[r] \ar[d]_{\simeq}& \Ell(X) \ar[d]^\simeq\\
K_1(\mathcal{M}\setminus X) \ar[r] & K_0(X) \ar[r] & K_0(\mathcal{M}).
}
$$
Since the map $K_0(X)\to K_0(\mathcal{M})$ preserves the
index (in $\mathbb{Z}$), we see that the index of $D$ on
$X$ is zero provided that $[D]\in \Ell(X)$ can be pulled back to
$\Ell_1(\mathcal{M}, X)$.

\begin{remark}
For nonsmooth $X$, the construction of this commutative
diagram is an open problem. (Even the exact sequence in $\Ell$-theory
is not known.)
\end{remark}

\section{Supplement. Proof of Theorem~\ref{luke2}}

Note that Theorem~\ref{luke2} is a strengthening of
the Luke theorem \cite{Luk1} for the index of operators with
operator-valued symbols. (The index theorem is obtained from
Eq.~\eqref{luuk} if $\partial X=\emptyset$ and we consider only the
indices of the corresponding $K$-homology elements.)

For the case of symbols homogeneous for large $|\xi|$,
the proof carries over from~\cite{Luk1}
word for word. The general case (nonhomogeneous symbols) is reduced
to the case of homogeneous symbols by the method suggested in~\cite{NSScS15},
where it was adapted to the index computation. The
reduction used in the present paper is based on the following standard fact
of Kasparov's $KK$-theory.
\begin{proposition}\label{chudo}
Let $P_t$ be a $*$-strongly continuous homotopy of bounded operators such that
$$
f[P_tP_t^*-1],\quad f[P_t^*P_t-1],\quad [f,P_t], \quad [f,P^*_t] \qquad \forall
f\in C_0(X^\circ)
$$
are norm continuous families of compact operators. Then the corresponding
element in $K$-homology does not depend on the parameter:
$$
[P_0]=[P_1]\in K_*(X^\circ).
$$
\end{proposition}
\begin{proof}
It follows from the conditions that the family $\{P_t\}$ defines an
operator in the $C_0(X^\circ)-C([0,1])$ bimodule $C_0(X^\circ\times
[0,1],L^2(H))$. This operator defines a homotopy (in the sense of
$KK$-theory) between $P_0$ and $P_1$ (e.g., see~\cite{Bla1}). Hence the invariance of the
$K$-homology element $[P_t]$ on $t$ follows from the equivalence of the
definitions of $K$-homology in terms of homotopies and in terms of operator homotopies.
\end{proof}

Let us apply this proposition. Without loss of generality, we can assume
that the symbol $p(x,\xi)$ is smooth up to the zero section in $T^*X$ and
normalized: $p^*(x,\xi)p(x,\xi)=1$ for large
$|\xi|$. Let $\psi(t)$, $t\ge0$, be a smooth positive function
such that
\begin{equation}\label{psi}
\psi(t)=\begin{cases}
1&\text{for } t<1,\\
1/t&\text{for } t>2.
\end{cases}
\end{equation}
Consider the family of symbols
\begin{equation}\label{homoe}
p_\e(x,\xi)=p(x,\xi\psi(\e|\xi|)).
\end{equation}
The following properties are obtained by a
straightforward computation
(cf.~\cite[Theorem~19.2.3]{Hor3}):
\begin{enumerate}
\item $p_0=p$. For $\e>0$, the symbol $p_\e$ is homogeneous for
large $|\xi|$,
and for $\e$ small it is elliptic.
\item $p_\e$ and $p^*_\e$ are uniformly bounded in the class of symbols
of compact variation in the fibers of $T^*X$
for $\e\in[0,1]$.
\item For small $\e>0$,
$p_\e p^*_\e-1$ and $p^*_\e p_\e-1$
are compactly supported and compact valued symbols independent of
$\e$.
\end{enumerate}

Consider the operator
\begin{equation}\label{operpeqe}
P_\e=p_\e\BL(x,-i\pd{}x\BR).
\end{equation}
(We fix an atlas of charts and a subordinate partition of
unity independent of $\e$ on $X$.) We claim that (for
small $\e\ge 0$)
\begin{enumerate}
\item[(a)] $\ind p_\e\in K_c(T^*X)$ is independent of $\e$.
\item[(b)]$P_\e$ satisfies the conditions of Proposition~\ref{chudo}.
\end{enumerate}
This will imply Theorem~\ref{luke2}, since, on the one hand, $p_\e$
is homogeneous at infinity for small $\e>0$, and
hence $[P_\e]=\varphi(\ind p_\e)$ (by the first part of the
proof) and on the other hand, the passage to the limit
as $\e\to0$ is possible by Proposition~\ref{chudo}.
Thus it remains to prove~(a) and~(b).

Claim (a) follows from the homotopy invariance of the index, since $p_\e$
changes under variations of $\e$ only in the complement of a large ball
$\{\abs{\xi}>R\}$, where $R\simeq1/\e$, and is invertible on this complement.

Claim (b) is proved as follows.

1) The families $P_\e$ and $P_\e^*$ are strongly continuous, since they are
uniformly bounded and each summand in their definitions in coordinate
patches on $X$ is strongly continuous on the set of functions whose
Fourier transform has compact support.

2) Given $f\in C_0(X^\circ)$, the operators $f(P_\e P^*_\e-1)$ and
$f(P^*_\e P_\e-1)$ are compact and continuously depend on $\e$. Indeed,
compactness is obvious, and the continuity
follows from the fact that their complete symbols and their derivatives in
local coordinates are uniformly continuous in $\e$ on compact subsets in
$\xi$, are uniformly bounded and decay to zero as $\xi\to\infty$, and
hence, are uniformly continuous in $\e$ for all $\xi$.

The compactness and continuity of the commutators
$[f,P_\e]$ and $[f,P^*_\e]$ can be proved along the same
lines.

The proof of Theorem~\ref{luke2} is
complete.

\end{document}